\documentclass[10pt,twoside,reqno]{amsart}
\usepackage{amssymb}
\usepackage{multirow}
\calclayout

\numberwithin{equation}{section}
\makeatletter

\renewcommand{\@secnumfont}{\bfseries}

\renewcommand{\section}{\@startsection{section}{1}%
  {0mm}{.7\linespacing\@plus\linespacing}{.5\linespacing}
  {\normalfont\bfseries\centering}}

\newcommand{\bibsection}{\@startsection{section}{1}%
  {0mm}{.7\linespacing\@plus\linespacing}{.5\linespacing}
  {\normalfont\scshape\centering}}

\renewcommand{\@biblabel}[1]{#1.}

\newtheorem{thm}{\bf Theorem}[section]

\begin{document}

\vspace{1.3cm}

\title {A note on degenerate Stirling polynomials of the second kind}

\author{Taekyun Kim}
\address{Department of Mathematics, Kwangwoon University, Seoul 139-701, Republic
	of Korea}
\email{tkkim@kw.ac.kr}


\markboth{\centerline{A note on degenerate Stirling polynomials of the second kind}}{\centerline{\scriptsize }}

\begin{abstract} 
In this paper, we consider the degenerate Stirling polynomials of the second kind which are derived from the generating function. In addition, we give some new identities for these polynomials.

\end{abstract}
\maketitle
\bigskip
\medskip
\section{\bf Introduction}
For $n \in \mathbb{N} \cup\{0\}$, as is well known, the Stirling number of the first kind is defined by
\begin{equation}\begin{split}\label{01}
(x)_0=1,\,\,(x)_n=x(x-1)\cdots(x-n+1)= \sum_{l=0}^n S_1(n,l) x^l,\,\,(n \geq 1).
\end{split}\end{equation}
Note that
\begin{equation}\begin{split}\label{02}
S_1(n+1,k)=S_1(n,k-1)-nS_1(n,k),\,\,(1 \leq k \leq n),\quad (\textnormal{see} \,\, [6]).
\end{split}\end{equation}
The Stirling number of the second kind is defined by 
\begin{equation}\begin{split}\label{03}
x^n = \sum_{l=0}^n S_2(n,l)(x)_l \,\,(n \geq 0),\quad (\textnormal{see} \,\, [1,2,4]).
\end{split}\end{equation}
From \eqref{03}, we note that
\begin{equation}\begin{split}\label{04}
S_2(n+1,k)= kS_2(n,k)+S_2(n,k-1),
\end{split}\end{equation}
where $1 \leq k \leq n$, (see [2,3,4,5]). The generating functions for $S_1(n,k)$ and $S_2(n,k)$, $(n,k \geq 0)$, are given by
\begin{equation}\begin{split}\label{05}
\frac{1}{k!} \big( \log(1+t)\big)^k = \sum_{n=k}^\infty S_1(n,k) \frac{t^n}{n!},
\end{split}\end{equation}
and 
\begin{equation}\begin{split}\label{06}
\frac{1}{k!}\big(e^t-1\big)^k = \sum_{n=k}^\infty S_2(n,k) \frac{t^n}{n!},\quad (\textnormal{see} \,\, [9]).
\end{split}\end{equation}
Now, we define the difference operator $\Delta$ as follows:
\begin{equation}\begin{split}\label{07}
\Delta f(x) = f(x+1) - f(x),\quad (\textnormal{see} \,\, [3,7]).
\end{split}\end{equation}
From \eqref{07}, we have
\begin{equation}\begin{split}\label{08}
\Delta^n f(x) = \sum_{k=0}^n {n \choose k} (-1)^{n-k} f(x+k),\,\,(n \in \mathbb{N} \cup \{0\}).
\end{split}\end{equation}
and
\begin{equation}\begin{split}\label{09}
f(x) \approx \sum_{k=0}^\infty {x \choose k} \Delta^k f(0),\quad (\textnormal{see} \,\, [3,7]).
\end{split}\end{equation}
From \eqref{08}, we note that
\begin{equation}\begin{split}\label{10}
\Delta^k 0^n = \begin{cases}S_2(n,k)&\text{if}\,\,n\geq k \\
0&\text{if}\,\,n <k.
\end{cases}
\end{split}\end{equation}
The Bernoulli polynomials are defined by the generating function to be
\begin{equation}\begin{split}\label{11}
\frac{t}{e^t-1}e^{xt} = \sum_{n=0}^\infty B_n(x) \frac{t^n}{n!},\quad (\textnormal{see} \,\, [ 1-10]).
\end{split}\end{equation}
When $x=0$, $B_n=B_n(0)$, $(n \geq 0)$, are called the Bernoulli numbers.

 The Euler polynomials are defined by the generating function as follows:
 \begin{equation}\begin{split}\label{12}
 \frac{2}{e^t+1}e^{xt} = \sum_{n=0}^\infty E_n(x) \frac{t^n}{n!},\quad (\textnormal{see} \,\, [ 1,2]).
 \end{split}\end{equation}
When $x=0$, $E_n=E_n(0)$ are Euler numbers. We observe that 
\begin{equation}\begin{split}\label{13}
\frac{2}{e^t+1} &= \left( \frac{e^t-1}{2} +1 \right)^{-1} = \sum_{l=0}^\infty \left( \frac{e^t-1}{2} \right)^l (-1)^l\\
&= \sum_{l=0}^\infty (-1)^l 2^{-l} l! \sum_{n=l}^\infty S_2(n,l) \frac{t^n}{n!}\\
&= \sum_{n=0}^\infty \left( \sum_{l=0}^n S_2(n,l) (-1)^l 2^{-l} l! \right) \frac{t^n}{n!}.
\end{split}\end{equation}
By \eqref{12} and \eqref{13}, we get
\begin{equation*}\begin{split}
E_n = \sum_{l=0}^n S_2(n,l) 2^{-l} l! (-1)^l,\,\,(n \geq 0).
\end{split}\end{equation*}
In [1], L. Carlitz consider the degenerate Bernoulli and Euler polynomials which are given by the generating function to be
\begin{equation}\begin{split}\label{14}
\frac{t}{(1+\lambda t)^{\frac{1}{\lambda }}-1}(1+\lambda t)^{\frac{x}{\lambda }} = \sum_{n=0}^\infty \beta_{n,\lambda }(x) \frac{t^n}{n!},
\end{split}\end{equation}
and
\begin{equation}\begin{split}\label{15}
\frac{2}{(1+\lambda t)^{\frac{1}{\lambda }}+1}(1+\lambda t)^{\frac{x}{\lambda }} = \sum_{n=0}^\infty \mathcal{E}_{n,\lambda }(x) \frac{t^n}{n!}.
\end{split}\end{equation}
Note that
\begin{equation*}\begin{split}
\lim_{\lambda \rightarrow 0} \beta_{n,\lambda}(x) =B_n(x),\,\, \lim_{\lambda \rightarrow 0} \mathcal{E}_{n,\lambda}(x) =E_n(x), \,\,(n \geq 0).
\end{split}\end{equation*}
In this paper, in the viewpoint \eqref{14} and \eqref{15}, we consider the degenerate Stirling polynomials of the second kind which are derived from the generating function. In addition, we give some new identities for these polynomials.

\section{Degenerate Stirling polynomials of the second kind}

Now, we define the Stirling polynomials of the second kind which are given by the generating function to be
\begin{equation}\begin{split}\label{16}
\frac{1}{k!} e^{xt} \big(e^t-1\big)^k = \sum_{n=k}^\infty S_2(n,k|x) \frac{t^n}{n!}.
\end{split}\end{equation}
Note that
\begin{equation}\begin{split}\label{17}
& \frac{1}{k!}e^{xt}\big(e^t-1\big)^k = \frac{1}{k!} \big(e^t-1\big)^k e^{xt} \\
&=\left( \sum_{l=k}^\infty S_2(l,k)\frac{t^l}{l!} \right) \left( \sum_{m=0}^\infty \frac{1}{m!} x^m t^m \right)\\
&= \sum_{n=k}^\infty \left\{ \sum_{l=k}^n {n \choose l} S_2(l,k) x^{n-l} \right\} \frac{t^n}{n!}
\end{split}\end{equation}
From \eqref{16} and \eqref{17}, we have
\begin{equation}\begin{split}\label{18}
S_2(n,k|x)= \sum_{l=k}^n {n \choose l} S_2(l,k) x^{n-l} ,\,\,(n,k \geq 0).
\end{split}\end{equation}
When $x=0$, we easily get $S_2(n,k|0)=S_2(n,k)$. Now, we consider the degenerate Stirling polynomials which are defined by the generating function as follows:
\begin{equation}\begin{split}\label{19}
\frac{1}{k!}(1+\lambda t)^{\frac{x}{\lambda }} \big( (1+\lambda t)^{\frac{1}{\lambda }}-1 \big)^k = \sum_{n=k}^\infty S_{2,\lambda }(n,k|x) \frac{t^n}{n!}.
\end{split}\end{equation}
Now, we observe that
\begin{equation}\begin{split}\label{20}
(1+\lambda t)^{\frac{x}{\lambda }} &= \sum_{l=0}^\infty {\frac{x}{\lambda } \choose l} \lambda ^l t^l = \sum_{l=0}^\infty \left(\frac{x}{\lambda }\right)_l \lambda ^l \frac{t^l}{l!}\\
&= \sum_{l=0}^\infty (x)_{l,\lambda } \frac{t^l}{l!}= \sum_{l=0}^\infty {x \choose l}_\lambda  t^l,
\end{split}\end{equation}
where
\begin{equation}\begin{split}\label{21}
(x)_{0,\lambda }=1,\,\,(x)_{l,\lambda }=x(x-\lambda )\cdots(x-(l-1)\lambda ),\,\,(l \geq 1),
\end{split}\end{equation}
and
\begin{equation}\begin{split}\label{22}
{x \choose l}_\lambda  = \frac{(x)_{l,\lambda }}{l!} = \frac{x(x-\lambda )\cdots(x-(l-1)\lambda )}{l!}.
\end{split}\end{equation}
For $n \in \mathbb{N}$, let us define $\lambda$-analogue of $n!$ as follows:
\begin{equation}\begin{split}\label{23}
(n)_\lambda ! &= n(n-\lambda )(n-2\lambda )\cdots(n-(n-1)\lambda )\\
&=(n)_{n,\lambda },
\end{split}\end{equation}
and
\begin{equation}\begin{split}\label{24}
{n \choose k}_\lambda  = \frac{(n)_\lambda !}{k!(n-k\lambda )_{n-k,\lambda }} = \frac{(n)_{k,\lambda }}{k!}, \,\,(n \geq k \geq 0).
\end{split}\end{equation}
From \eqref{19}, we have
\begin{equation}\begin{split}\label{25}
&\frac{1}{k!}(1+\lambda t)^{\frac{x}{\lambda }} \big( (1+\lambda t)^{\frac{1}{\lambda }}-1 \big)^k \\
&=\frac{1}{k!}(1+\lambda t)^{\frac{x}{\lambda }}  \sum_{l=0}^k {k \choose l}(-1)^{k-l} (1+\lambda t)^{\frac{l}{\lambda }}\\
&= \frac{1}{k!} \sum_{l=0}^k {k \choose l} (-1)^{k-l} (1+\lambda t)^{\frac{l+x}{\lambda }}\\
&= \sum_{n=0}^\infty \left( \frac{n!}{k!} \sum_{l=0}^k  {k \choose l}(-1)^{k-l} {l+x \choose n}_\lambda  \right) \frac{t^n}{n!}.
\end{split}\end{equation}
Therefore, by \eqref{19} and \eqref{25}, we obtain the following theorem.

\begin{thm}
For $n,k \geq 0$, we have
\begin{equation*}\begin{split}
\frac{n!}{k!} \sum_{l=0}^k  {k \choose l}(-1)^{k-l} {l+x \choose n}_\lambda = \begin{cases}
S_{2,\lambda }(n,k|x),&\text{if}\,\,n \geq k
\\0,&\text{if}\,\,n<k.
\end{cases}
\end{split}\end{equation*}
\end{thm}
From \eqref{19}, we have
\begin{equation}\begin{split}\label{26}
&\frac{1}{k!}(1+\lambda t)^{\frac{x}{\lambda }} \big( (1+\lambda t)^{\frac{1}{\lambda }}-1 \big)^k \\
&=\frac{1}{k!} e^{\frac{x}{\lambda }\log(1+\lambda t)}\left(e^{\frac{1}{\lambda }\log(1+\lambda t)}-1\right)^k\\
&= \frac{1}{k!}\sum_{l=0}^k {k \choose l} (-1)^{k-l} e^{\frac{x+l}{\lambda }\log(1+\lambda t)}\\
&= \frac{1}{k!}\sum_{l=0}^k {k \choose l} (-1)^{k-l} \sum_{m=0}^\infty \left( \frac{x+l}{\lambda }\right)^m \frac{1}{m!} \big(\log(1+\lambda t)\big)^m\\
&= \frac{1}{k!}\sum_{l=0}^k {k \choose l} (-1)^{k-l} \sum_{m=0}^\infty \frac{(x+l)^m}{\lambda ^m} \sum_{n=m}^\infty S_1(n,m) \frac{\lambda ^n t^n}{n!}\\
&= \sum_{n=0}^\infty \frac{1}{k!} \sum_{l=0}^k {k \choose l} (-1)^{k-l} \sum_{m=0}^n \lambda ^{n-m} S_1(n,m) (x+l)^m \frac{t^n}{n!}\\
&= \sum_{n=0}^\infty \left\{ \sum_{m=0}^n \left( \frac{1}{k!} \sum_{l=0}^k {k \choose l} (-1)^{k-l} (x+l)^m \right) \lambda ^{n-m} S_1(n,m) \right\} \frac{t^n}{n!}\\
&= \sum_{n=0}^\infty \left( \sum_{m=0}^n \frac{1}{k!} \Delta^k x^m \lambda ^{n-m} S_1(n,m) \right) \frac{t^n}{n!}.
\end{split}\end{equation}
Therefore, by \eqref{19} and \eqref{26}, we obtain the following theorem.
\begin{thm} For $n,k \geq 0$, we have
\begin{equation*}\begin{split}
\sum_{m=0}^n \frac{1}{k!} \Delta^k x^m \lambda ^{n-m} S_1(n,m)\begin{cases}
S_{2,\lambda }(n,k|x),&\text{if}\,\,n \geq k
\\0,&\text{if}\,\,n<k.
\end{cases}
\end{split}\end{equation*}
\end{thm}
By \eqref{08}, we easily get
\begin{equation}\begin{split}\label{27}
\Delta^k x^{m+1} &= \sum_{l=0}^k {k \choose l}(x+l)^{m+1} (-1)^{k-l} = \sum_{l=0}^k {k \choose l} (-1)^{k-l} (x+l)^m (x+l) \\
&=x \sum_{l=0}^k {k \choose l} (x+l)^m (-1)^{k-l} + \sum_{l=0}^k {k \choose l} (x+l)^m l (-1)^{k-l}\\
&= x \Delta^k x^m + k \sum_{l=1}^k {k-1 \choose l-1} (x+l)^m (-1)^{k-l}\\
&= x \Delta^k x^m + k \sum_{l=1}^k \left\{ {k \choose l} - {k-1 \choose l} \right\} (x+l)^m (-1)^{k-l}\\
&= x \Delta^k x^m + k \sum_{l=0}^k \left\{ {k \choose l} - {k-1 \choose l} \right\} (x+l)^m (-1)^{k-l}\\
&= x \Delta^k x^m + k \big( \Delta^k x^m + \Delta^{k-1} x^m \big).
\end{split}\end{equation}
From \eqref{02} and \eqref{27}, we have
\begin{equation}\begin{split}\label{28}
S_{2,\lambda }(n+1,k|x) &= \sum_{m=0}^{n+1} \frac{1}{k!} \Delta^k x^m \lambda ^{n+1-m} S_1(n+1,m)\\
 &= \sum_{m=1}^{n+1} \frac{1}{k!} \Delta^k x^m \lambda ^{n+1-m} S_1(n+1,m)\\
 &= \sum_{m=1}^{n+1} \frac{1}{k!} \Delta^k x^m \lambda ^{n+1-m} \big( S_1(n,m-1) - n S_1(n,m)\big) \\
 &= \sum_{m=1}^{n+1} \frac{1}{k!} \Delta^k x^m \lambda ^{n+1-m} S_1(n,m-1)\\
 &\,\, - n\lambda  \sum_{m=0}^n \frac{1}{k!} \Delta^k x^m \lambda ^{n-m} S_1(n,m)\\
 &= \sum_{m=1}^{n+1} \frac{1}{k!} \Delta^kx^m \lambda ^{n+1-m} S_1(n,m-1) - n\lambda S_{2,\lambda }(n,k|x).
\end{split}\end{equation}
Now, we observe that
\begin{equation}\begin{split}\label{29}
&\sum_{m=1}^{n+1} \frac{1}{k!} \Delta^k x^m \lambda^{n+1-m} S_1(n,m-1) = \sum_{m=0}^n \frac{1}{k!} \Delta^k x^{m+1} \lambda ^{n-m} S_1(n,m)\\
&= \sum_{m=0}^n \frac{1}{k!} \left\{ x \Delta^k x^m + k\big( \Delta^k x^m + \Delta^{k-1} x^m \big) \right\} \lambda ^{n-m} S_1(n,m)\\
&= x \frac{1}{k!} \sum_{m=0}^n \Delta^k x^m \lambda ^{n-m} S_1(n,m) + k \frac{1}{k!} \Delta^k x^m \lambda ^{n-m} S_1(n,m) \\
&\,\,+ \frac{1}{(k-1)!} \sum_{m=0}^n \Delta^{k-1} x^m \lambda ^{n-m} S_1(n,m)\\
&= (x+k) S_{2,\lambda }(n,k|x) + S_{2,\lambda }(n,k-1|x),\,\,(1 \leq k \leq n).
\end{split}\end{equation}
Therefore, by \eqref{28} and \eqref{29}, we obtain the following theorem.
\begin{thm}
For $1 \leq k \leq n$, we have
\begin{equation*}\begin{split}
S_{2,\lambda }(n+1,k|x) = (x+k)S_{2,\lambda }(n,k|x)+S_{2,\lambda }(n,k-1|x)-n\lambda S_{2,\lambda }(n,k|x).	
	\end{split}\end{equation*}	
\end{thm}
Note that
\begin{equation}\begin{split}\label{30}
\lim_{\lambda  \rightarrow 0} S_{2,\lambda }(n+1,k|x)  = (x+k) S_2(n,k|x) + S_2(n,k-1|x).
\end{split}\end{equation}
Thus, by \eqref{30} we get
\begin{equation}\begin{split}\label{31}
S_2(n+1,k|x) = (x+k) S_2(n,k|x) + S_2(n,k-1|x),
\end{split}\end{equation}
where $1 \leq k \leq n$.
As is known, the higher-order Carlitz degenerate Euler polynomials are defined by the generating function as follows:
\begin{equation}\begin{split}\label{32}
\left( \frac{2}{(1+\lambda t)^{\frac{1}{\lambda }}+1} \right)^r (1+\lambda t)^{\frac{x}{\lambda }} = \sum_{n=0}^\infty \mathcal{E}_{n,\lambda }^{(r)}(x) \frac{t^n}{n!},\,\,(r \in \mathbb{N}).
\end{split}\end{equation}
From \eqref{32}, we note that
\begin{equation}\begin{split}\label{33}
&\left( \frac{2}{(1+\lambda t)^{\frac{1}{\lambda }}+1} \right)^r (1+\lambda t)^{\frac{x}{\lambda }} =2^r \big( (1+\lambda t)^{\frac{1}{\lambda }}+1\big)^{-r} (1+\lambda t)^{\frac{x}{\lambda }}\\&= \left( \frac{(1+\lambda t)^{\frac{1}{\lambda }}-1}{2}+1\right)^{-r}  (1+\lambda t)^{\frac{x}{\lambda }}\\
&=\sum_{l=0}^\infty {r+l-1 \choose l} 2^{-l} (-1)^l \big( (1+\lambda t)^{\frac{1}{\lambda }}-1\big)^l (1+\lambda t)^{\frac{x}{\lambda }}\\
&=\sum_{l=0}^\infty {r+l-1 \choose l} 2^{-l} (-1)^l l! \sum_{n=l}^\infty S_{2,\lambda }(n,l|x) \frac{t^n}{n!}\\
&= \sum_{n=0}^\infty \left( \sum_{l=0}^n {r+l-1 \choose l}2^{-l}(-1)^l l! S_{2,\lambda }(n,l|x) \right) \frac{t^n}{n!}.
\end{split}\end{equation}
Therefore, by \eqref{32} and \eqref{33}, we obtain the following theorem.
\begin{thm}
For $n \geq 0$, we have
\begin{equation*}\begin{split}
\mathcal{E}_{n,\lambda }^{(r)}(x) = \sum_{l=0}^n {r+l-1 \choose l}2^{-l}(-1)^l l! S_{2,\lambda }(n,l|x).
\end{split}\end{equation*}
\end{thm}
For $r \in \mathbb{N}$, the $r$-Whitney numbers of the second kind $W_{m,r}(n,k)$ are defined by
\begin{equation}\begin{split}\label{34}
(mx+r)^n = \sum_{k=0}^n m^k W_{m,r}(n,k) (x)_k,\quad (\textnormal{see} \,\, [8]).
\end{split}\end{equation}
where $n,r \in \mathbb{N}\cup \{0\}$ and $m \in \mathbb{N}$.
The generating function of Whitney numbers is given by 
\begin{equation}\begin{split}\label{35}
\frac{1}{m^k k!} e^{rt} (e^{mt}-1)^k = \sum_{n=k}^\infty W_{m,r}(n,k) \frac{t^n}{n!}.
\end{split}\end{equation}
From \eqref{35}, we note that
\begin{equation}\begin{split}\label{36}
\frac{1}{m^k k!} e^{rt} (e^{mt}-1)^k &= \frac{1}{m^k} \left( \sum_{l=0}^\infty r^l \frac{t^l}{l!}\right) \left( \sum_{i=k}^\infty S_2(i,k) \frac{m^i t^i}{i!} \right)\\
&= \frac{1}{m^k} \sum_{n=k}^\infty \left( \sum_{i=k}^n S_2(i,k) {n \choose i} r^{n-i} m^i \right) \frac{t^n}{n!}.
\end{split}\end{equation}
Thus, by \eqref{35} and \eqref{36}, we get
\begin{equation}\begin{split}\label{37}
W_{m,r}(n,k) = \sum_{i=k}^n  {n \choose i} r^{n-i} S_2(i,k)  m^{i-k}.
\end{split}\end{equation}
By \eqref{34}, we easily get 
\begin{equation}\begin{split}\label{38}
(mx+r)^n = \sum_{k=0}^n\left( \sum_{i=k}^n {n \choose i} r^{n-i} m^i S_2(i,k) \right) (x)_k,
\end{split}\end{equation}
where $n,r \in \mathbb{N} \cup \{0\}$ and $m \in \mathbb{N}$.

Note that $W_{1,0}(n,k) = S_2(n,k)$ and $W_{1,r}(n,k) = S_2(n,k|r)$. From \eqref{35}, we note that

\begin{equation}\begin{split}\label{39}
&\frac{1}{m^k k!} e^{rt} \big( e^{mt}-1\big)^k = \frac{1}{m^k k!} e^{rt} \sum_{l=0}^k {k \choose l} (-1)^{k-l} e^{lm}\\
&=  \frac{1}{m^k k!}  \sum_{l=0}^k {k \choose l} (-1)^{k-l} e^{(lm+r)t} = \frac{1}{m^k k!} \sum_{l=0}^k {k \choose l} (-1)^{k-l} \sum_{n=0}^\infty (lm+r)^n \frac{t^n}{n!}\\
&= \sum_{n=0}^\infty \left( \frac{1}{m^k k!} \sum_{l=0}^k {k \choose l} (-1)^{k-l} (lm+r)^n \right) \frac{t^n}{n!}\\
&= \sum_{n=0}^\infty m^{n-k} \left( \frac{1}{k!} \sum_{l=0}^k {k \choose l} (-1)^{k-l} \big( l + \tfrac{r}{m} \big)^n \right) \frac{t^n}{n!}\\
&= \sum_{n=0}^\infty \left( m^{n-k} \frac{1}{k!} \Delta^k \left( \frac{r}{m} \right)^n \right) \frac{t^n}{n!}.
\end{split}\end{equation}
Therefore, by \eqref{35} and \eqref{39}, we obtain the following theorem.
\begin{thm}
For $m \in \mathbb{N}$ and $n,k,r \in \mathbb{N} \cup \{0\}$, we have
\begin{equation*}\begin{split}
W_{m,r}(n,k) = \begin{cases}
 m^{n-k} \frac{1}{k!} \Delta^k \left( \frac{r}{m} \right)^n,&\text{if}\,\,n \geq k,\\
 0,&\text{if}\,\,n<k.
 \end{cases}
\end{split}\end{equation*}
\end{thm}
For $1 \leq k \leq n$, by Theorem 5, we easily get
\begin{equation}\begin{split}\label{40}
W_{m,r}(n+1,k) = (r+mk)W_{m,r}(n,k) + W_{m,r}(n,k-1).
\end{split}\end{equation}
We consider the degenerate Whitney numbers which are defined by the generating function to be
\begin{equation}\begin{split}\label{41}
\frac{1}{m^k k!} (1+\lambda t)^{\frac{r}{\lambda }} \big( (1+\lambda t)^{\frac{m}{\lambda }}-1 \big)^k = \sum_{n=k}^\infty W_{m,r}(n,k|\lambda ) \frac{t^n}{n!},
\end{split}\end{equation}
where $m \in \mathbb{N}$, $n,k,r \in \mathbb{N} \cup \{0\}$. Note that $\lim_{\lambda  \rightarrow 0} W_{m,r}(n,k|\lambda ) = W_{m,r}(n,k).$ Now, we observe that
\begin{equation}\begin{split}\label{42}
&\frac{1}{m^k k!} (1+\lambda t)^{\frac{r}{\lambda }}\big((1+\lambda t)^{\frac{m}{\lambda }}-1\big)^k\\
&= \frac{1}{m^k k!} \sum_{l=0}^k {k \choose l} (-1)^{k-l} (1+\lambda t)^{\frac{ml}{\lambda }} (1+\lambda t)^{\frac{r}{\lambda }}\\
&= \sum_{n=0}^\infty \left( \frac{n!}{k!} \sum_{l=0}^k {k \choose l} (-1)^{k-l} {ml+r \choose n}_\lambda  \right) \frac{t^n}{n!}.
\end{split}\end{equation}
Therefore, by \eqref{41} and \eqref{42}, we obtain the following theorem.
\begin{thm}
For $m \in \mathbb{N}$, and $n,k,r \in \mathbb{N} \cup \{0\}$, we have
\begin{equation*}\begin{split}
 \frac{n!}{k!} \sum_{l=0}^k {k \choose l} {ml+r \choose n}_\lambda  (-1)^{k-l}= \begin{cases}
W_{m,r}(n,k|\lambda ),&\text{if}\,\,n \geq k,\\
 0,&\text{if}\,\,n<k.
 \end{cases}
\end{split}\end{equation*}
\end{thm}
From \eqref{41}, we note that
\begin{equation}\begin{split}\label{43}
&\frac{1}{m^k k!} (1+\lambda t)^{\frac{r}{\lambda }} \big( (1+\lambda t)^{\frac{m}{\lambda }}-1\big)^k = \frac{1}{m^k k!} e^{\frac{r}{\lambda }\log(1+\lambda t)} \big( e^{\frac{m}{\lambda }\log(1+\lambda t)}-1\big)^k\\
&= \frac{1}{k! m^k} \sum_{l=0}^k {k \choose l} (-1)^{k-l} e^{\frac{lm+r}{\lambda }\log(1+\lambda t)}\\
&= \frac{1}{k! m^k} \sum_{l=0}^k {k \choose l} (-1)^{k-l} \sum_{j=0}^\infty  \left( \frac{lm+r}{\lambda } \right)^j \frac{1}{j!} \big( \log(1+\lambda t)\big)^j\\
&= \frac{1}{k! m^k} \sum_{l=0}^k {k \choose l} (-1)^{k-l} \sum_{j=0}^\infty \left( \frac{lm+r}{\lambda }\right)^j \sum_{n=j}^\infty S_1(n,j) \frac{\lambda ^n t^n}{n!}\\
&=\frac{1}{k! m^k} \sum_{n=0}^\infty \left( \sum_{j=0}^n \lambda ^{n-j} S_1(n,j) \sum_{l=0}^k {k \choose l} (-1)^{k-l} (lm+r)^j \right) \frac{t^n}{n!}\\
&=\frac{1}{k!m^k} \sum_{n=0}^\infty \left( \sum_{j=0}^n \lambda ^{n-j} m^j S_1(n,j) \sum_{l=0}^k {k \choose l} (-1)^{k-l} \big(l+\tfrac{r}{m}\big)^j\right) \frac{t^n}{n!}\\
&= \sum_{n=0}^\infty \left( \sum_{j=0}^n \lambda ^{n-j} m^j S_1(n,j) \frac{1}{k! m^k} \Delta^k \left(\frac{r}{m}\right)^j \right) \frac{t^n}{n!}.
\end{split}\end{equation}
Therefore, by \eqref{04} and \eqref{43}, we obtain the following theorem.
\begin{thm}
For $n \in \mathbb{N}$, and $n,k,r \in \mathbb{N} \cup \{0\}$, we have
\begin{equation*}\begin{split}
 \frac{1}{k! m^k}\sum_{j=0}^n \lambda ^{n-j}  S_1(n,j) m^j \Delta^k \left(\frac{r}{m}\right)^j=
 \begin{cases}
W_{m,r}(n,k|\lambda ),&\text{if}\,\,n \geq k,\\
 0,&\text{if}\,\,n<k.
 \end{cases}
\end{split}\end{equation*}
\end{thm}
For $1\leq k \leq n$, by Theorem 7, we get
\begin{equation}\begin{split}\label{44}
&W_{m,r}(n+1,k|\lambda )= \frac{1}{k! m^k}\sum_{j=0}^{n+1} \lambda ^{n+1-j}  S_1(n+1,j) m^j \Delta^k \left(\frac{r}{m}\right)^j\\
&= \frac{1}{k! m^k}\sum_{j=1}^{n+1} \lambda ^{n+1-j} \Big\{ S_1(n,j-1)-nS_1(n,j)\Big\}  m^j \Delta^k \left(\frac{r}{m}\right)^j\\
&=  \frac{1}{k! m^k}\sum_{j=1}^{n+1} \lambda ^{n+1-j}S_1(n,j-1)m^j \Delta^k \left(\frac{r}{m}\right)^j - \frac{n\lambda }{k! m^k}\sum_{j=1}^{n+1} \lambda ^{n-j}  S_1(n,j) m^j \Delta^k \left(\frac{r}{m}\right)^j\\
&=  \frac{1}{k! m^k}\sum_{j=0}^{n} \lambda ^{n-j}S_1(n,j-1)m^{j+1} \Delta^k \left(\frac{r}{m}\right)^{j+1} \\
&\,\,-  \frac{n\lambda }{k! m^k}\sum_{j=0}^{n} \lambda ^{n-j}  S_1(n,j) m^j \Delta^k \left(\frac{r}{m}\right)^j\\
&= \frac{m}{k! m^k} \sum_{j=0}^n \lambda^{n-j} S_1(n,j) m^j \Delta^k \left(\frac{r}{m}\right)^{j+1} - n\lambda W_{m,r}(n,k|\lambda).
\end{split}\end{equation}
Now, we observe that

\begin{equation*}\begin{split}
&\Delta^k \left(\frac{r}{m}\right)^{j+1} = \sum_{l=0}^k {k \choose l} (-1)^{k-l} \Big(l+\frac{r}{m}\Big)^{j+1} = \sum_{l=0}^k {k \choose l} (-1)^{k-l}  \Big(l+\frac{r}{m}\Big)^j  \Big(l+\frac{r}{m}\Big) \\
&=\sum_{l=0}^k {k \choose l}(-1)^{k-l} \Big(l+\frac{r}{m}\Big) ^j l + \frac{r}{m} \sum_{l=0}^k {k \choose l} (-1)^{k-l}  \Big(l+\frac{r}{m}\Big) ^j\\
&= \frac{r}{m} \Delta^k \left( \frac{r}{m} \right)^j + k \sum_{l=1}^k {k-1 \choose l-1} (-1)^{k-l} \Big(l+\frac{r}{m}\Big) ^j\\
&= \frac{r}{m} \Delta^k \left( \frac{r}{m} \right)^j + k \sum_{l=0}^k \left\{ {k \choose l} - {k-1 \choose l}\right\} (-1)^{k-l} \Big(l+\frac{r}{m}\Big) ^j\\
\end{split}\end{equation*}
\begin{equation}\begin{split}\label{45}
&=\frac{r}{m} \Delta^k \left( \frac{r}{m} \right)^j + k \sum_{l=0}^k {k \choose l} (-1)^{k-l} \Big(l+\frac{r}{m}\Big) ^j + k \sum_{l=0}^{k-1} {k -1\choose l} (-1)^{k-1-l} \Big(l+\frac{r}{m}\Big) ^j\\
&=\frac{r}{m} \Delta^k \left( \frac{r}{m} \right)^j + k \Delta^k \left( \frac{r}{m} \right)^j+ k \Delta^{k-1} \left( \frac{r}{m} \right)^j.
\end{split}\end{equation}
From \eqref{44} and \eqref{45}, we note that
\begin{equation*}\begin{split}
&W_{m,r}(n+1,k|\lambda ) = \frac{m}{k! m^k} \sum_{j=0}^m m^j S_1(n,j) \lambda ^{n-j} \Delta^k \left( \frac{r}{m} \right)^{j+1} - n\lambda W_{m,r}(n,k|\lambda) \\
&= \frac{m}{k! m^k} \sum_{j=0}^m m^j S_1(n,j) \lambda ^{n-j}\left\{ \frac{r}{m} \Delta^m
\left( \frac{r}{m} \right)^j + k \Delta^k \left( \frac{r}{m} \right)^j + k \Delta^{k-1} \left( \frac{r}{m} \right)^j \right\}\\
&\,\,  - n\lambda W_{m,r}(n,k|\lambda) \\
&=\frac{r}{k! m^k} \sum_{j=0}^m m^j S_1(n,j) \lambda^{n-j} \Delta^m \left( \frac{r}{m} \right)^j + \frac{mk}{k! m^k} \sum_{j=0}^m m^j S_1(n,j) \lambda^{n-j} \Delta^k \left( \frac{r}{m} \right)^j\\
&\,\,+ \frac{1}{(k-1)!m^{k-1}} \sum_{j=0}^m m^j S_1(n,j) \lambda^{n-j} \Delta^{k-1} \left( \frac{r}{m} \right)^j - n \lambda W_{m,r}(n,k|\lambda)\\
&= r W_{m,r}(n,k|\lambda) +mkW_{m,r}(n,k|\lambda) +W_{m,r}(n,k-1|\lambda) -n \lambda W_{m,r}(n,k|\lambda) .
\end{split}\end{equation*}
Therefore, we obtain the following theorem.
\begin{thm}
For $1 \leq k \leq n$, we have
\begin{equation*}\begin{split}
W_{m,r}(n+1,k|\lambda ) = (r+mk)W_{m,r}(n,k|\lambda) + W_{m,r}(n,k-1|\lambda)- n\lambda W_{m,r}(n,k|\lambda) .
\end{split}\end{equation*}
\end{thm}

\noindent {\bf{Remark.}} From \eqref{20}, we note that
\begin{equation}\begin{split}\label{46}
(1+ \lambda t)^{\frac{x+y}{\lambda }}&= (1+\lambda t)^{\frac{x}{\lambda }} (1+\lambda t)^{\frac{y}{\lambda }} = \left( \sum_{l=0}^\infty {x \choose l}_\lambda  t^l \right) \left( \sum_{m=0}^\infty {y \choose m}_\lambda  t^m \right)\\
&= \sum_{n=0}^\infty \left( \sum_{m=0}^n {y \choose m}_\lambda  {x \choose n-m}_\lambda  \right) t^n.
\end{split}\end{equation}
By \eqref{46}, we easily get
\begin{equation}\begin{split}\label{47}
(1+ \lambda t)^{\frac{x+y}{\lambda }} = \sum_{n=0}^\infty {x+y \choose n}_\lambda t^n.
\end{split}\end{equation}
Comparing the coefficients on the both sides of \eqref{46} and \eqref{47}, we have
\begin{equation*}\begin{split}
\sum_{m=0}^n {y \choose m}_\lambda  {x \choose n-m}_\lambda  = {x+y \choose n}_\lambda,\,\,(n \geq 0).
\end{split}\end{equation*}

\end{document}